[1994/12/01]
\documentclass[12pt]{article}
\title{A row space method for solving a system of linear equations}
\author{Michael F. Zimmer}
\date{\today}	

\usepackage{amsmath,amsthm,url}

\begin{document}
\maketitle

\section*{Abstract}

A new algorithm is presented for computing a direct solution
to a system of consistent linear equations.  It produces a 
minimum norm particular solution, a generalized inverse (of type $\{124\}$), 
and a null space projection operator.  In addition, the 
algorithm permits an online formulation so that computations may 
proceed as the data become available.  The algorithm does 
not require the solution of a triangular system of equations, 
nor does it rely on block partitioned matrices.\\

\noindent \emph{AMS classification:} 15A06, 15A09\\

\noindent \emph{Key Words:}  row space; linear equations; system of linear equations; generalized inverse; minimum norm; null space projection;
direct solution; online solution

\section{Introduction}

The problem of solving a system of linear equations
is widespread across mathematics, science and engineering.
When there are $m$ equations in $n$ unknowns, the equations are
written as $Ax=b$, where $A=(a_{ij})$ is the coefficient matrix, 
$x=(x_j)$ is the vector of unknowns, and $b=(b_i)$.  (The indices
are $i=1, \dots, m$ and $j=1, \dots ,n$.)
The entries in $A$, $x$ and $b$ are in a field ${\cal F}$, which is assumed to be endowed with
an inner product.  The examples in this paper assume ${\cal F}={\cal C}$.

It is known that a solution to $Ax=b$ may be written as $x = x_p + x_h$,
where the particular solution $x_p$ satisfies $Ax_p=b$
and the homogeneous solution $x_h$ satisfies $Ax_h=0$.
Note that $x_h$ is an element of the null space of $A$ (i.e., ${\cal N}(A)$).
When the dimension of ${\cal N}(A)$ is greater than zero
or when $A$ is rectangular, the inverse
to $A$ does not exist, and hence it's no longer possible to
write the solution in the cogent form $x = A^{-1}b$.
Nevertheless, when $A^{-1}$ doesn't exist, one may still use
the generalized inverse ($G$), which allows the two parts
of the solution to be written as:
\begin{align*}
x_p & = Gb \\
x_h & = Py \\
P & = I_n - GA
\end{align*}
where $I_n$ is an n-by-n identity matrix, $P$ is a null space
projection operator and $y$ is an arbitrary n-dimensional vector.  
(By definition, $Py \in {\cal N}(A)$ for any $y \in {\cal F}^n$.)
Different properties of the solution $Gb$ can be inferred by the
type of the generalized inverse.  (This is discussed in the literature \cite{penrose1,penrose2,campbell,benis,bapat}
and will not be reviewed here.)  The $G$ that will be computed
here is of type $\{124\}$, which among other things yields a solution $Gb$
which has minimum Euclidean norm.

There are other interesting features, such as the fact that
the algorithm may be formulated in an \emph{online mode},
which is defined to mean that the solution $x$ can be
formed while the input data ($A$ and $b$) are still being acquired.
This is useful in cases when it takes a relatively long 
time to access/create the input data.

\section{The algorithm}

The solution will require $A$ to have its rows orthonormalized,
perhaps using the usual Gram-Schmidt orthogonalization procedure.
However, before discussing a modified orthonormalization procedure, 
it is helpful to define the following.\\

\noindent{\it Definition}: A {\bf quasi-orthonormal} list of vectors
consists of vectors with norm equal to 1 or 0.
Those vectors of norm equal to 1 are mutually orthogonal.\\ 

\noindent{\it Definition}: The {\bf row orthonormalization procedure} 
(ROP) of a matrix $Q$ is
performed by left multplying it by a series of matrices $M_s$ ($s=1,2,3,...$)
which affect an orthonormalization procedure on the nonzero rows of $Q$.
The result is that the rows of $Q$ form a quasi-orthonormal list
of vectors.\\ 

Thus, after applying the ROP, the rows are either zero
or part of an orthonormal set.  
(Pseudocode for the ROP may be found in Appendix A.)
The ROP will here be applied to $A$ as it exists in the context
of the equation $Ax=b$, causing the same operations to be applied to $b$ as well.
It is implemented with the nonsingular, m-by-m
matrices $M_s$ which results in
\begin{equation}
(\cdots M_2 M_1) A x = (\cdots M_2 M_1)b
\end{equation}
The $M_i$ are applied to $A$, causing it to become 
$A' = MA$, where $M = (\cdots M_2 M_1 )$.  (Note that it is not
required that the $M_s$ represent elementary row operations.)
There are two versions for how to view the right-hand side of the previous equation: 
(1) apply the $M_s$ to $b$ so that it becomes $b'$ (where $b'$ equals $Mb$);
(2) accumulate the $M_s$ as a prefactor so that it becomes $Mb$.
To summarize, after applying the ROP to
$Ax=b$, these two variations may be written as 
\begin{equation}\label{var1}
A' x = b'
\end{equation}
or
\begin{equation}\label{var2}
A' x = M b
\end{equation}
The significance of the difference has mainly to do with computer
implementation issues, with respect to storage and reusability.\\

\noindent{\it Definition}:
Let $W_Q$ be a set which consists of the indices of the non-zero 
rows of an m-by-m matrix $Q$.  Note that $W_Q \subseteq \{1,2,...,m\}$.
Define an {\bf index matrix} by the following:
\[({\cal I}_Q)_{ij}= \begin{cases}
1 & \text{if $i=j$ and $i \in W_Q$}, \\
0 & \text{otherwise.}
\end{cases} \]

This matrix is essentially just an identity matrix, 
except now the i-th row has a 1 only if $i \in W_Q$.
This matrix is immediately applicable to $A'$, whose rows
comprise a quasi-orthonormal list of vectors.  The following identities
may be easily verified:
\begin{align}\label{AAstar}
A'~ (A')^*   & =  {\cal I}_{A'} \\
A'           & =  {\cal I}_{A'} ~ A' \\
(A')^*       & =  (A')^* ~ {\cal I}_{A'} \\
{\cal I}_{A'} ~ b'  & =  b'
\end{align}
where the superscript * represents a conjugate transpose.\\

\noindent{\it Lemma}:
An arbitrary vector $x \in {\cal F}^n$ may be expressed as
\begin{equation}\label{paramsoln}
x = (A')^*w + x_h
\end{equation}
where $A$ is an m-by-n matrix, $w \in {\cal F}^m$ and $x_h$ satisfies $Ax_h = 0$.\\

\begin{proof}
Beginning with $A'x = b'$, observe that $A': {\cal F}^n \rightarrow {\cal F}^m$
is a linear transformation.  It is then known (Thm. 18.3 \cite{bowen})
\begin{equation}
{\cal F}^n = \text{range}[ (A')^*] \oplus \text{null}[A']
\end{equation}
This direct sum decomposition may be used to rewrite an $x \in {\cal F}^n$ as
\begin{equation}\label{xdecomp}
x = w_1 + w_2
\end{equation}
where 
\begin{align}
w_1 & \in \text{range}[ (A')^*] \\
w_2 & \in \text{null}[A']
\end{align}
However, every vector that is in the range of some matrix $Q$ may be expressed as $Qv$,
for an appropriate $v$.  Thus, for some $w \in {\cal F}^m$, it follows that $w_1$
may be written as 
\begin{equation}\label{w1}
w_1 = (A')^* w
\end{equation}
Also, because $M$ is nonsingular, $\text{null}[A'] = \text{null}[A]$,
so that $w_2 \in \text{null}[A]$.
Substituting \eqref{w1} into \eqref{xdecomp} gives the desired result.
\end{proof}
What remains is to determine $w$ and to find a means of computing $x_h$.\\

\noindent{\it Theorem}:
A solution to $Ax=b$ is 
\begin{equation}\label{mainsoln}
x = (A')^*b' + x_h
\end{equation}
where all variables are as defined earlier.\\
\begin{proof}
The proposed solution is verified by substitution into $Ax=b$
\begin{align*}
Ax & = A(A')^* b' + Ax_h\\
   & = M^{-1} A'(A')^* b' \\
   & = M^{-1} {\cal I}_{A'} b' \\
   & = M^{-1} b' \\
   & = b
\end{align*}
\end{proof}

One could also take a constructive approach toward obtaining
the above solution.  Upon substituting \eqref{paramsoln} into $A'x=b'$
one obtains
\begin{equation*}
A'(A')^*w = b'.
\end{equation*}
which upon using \eqref{AAstar} becomes ${\cal I}_{A'} ~w = b'$. 
Equation \eqref{paramsoln} may now be re-expressed as 
\begin{align*}
x & = (A')^* ~w + x_h\\
  & = (A')^* {\cal I}_{A'} ~w + x_h\\
  & = (A')^* ~b' + x_h
\end{align*}

\subsection*{Discussion}

One of the first things to notice is that $(A')^* b'$
is a particular solution to $Ax=b$, as well as a
minimum norm solution.  The reason for the latter is because
it is an element of $\text{range}[ (A')^*]$ which is the orthogonal complement
of the null space of $A'$.  In other words, the particular solution
is orthogonal to the null space.  The minimum norm nature of the particular
solution will be revisited when the Penrose identities are checked.

The first variation (Eqn. \eqref{var1}) transforms $Ax=b$ into $A'x=b'$.
In terms of augmented matrices, one applies the ROP to $[A|b]$
to produce $[A'|b']$, where the rows of $A'$ form a quasi-orthonormal
list of vectors.  The solution is formed as
\begin{align*}
x_p & = (A')^* b' \\
x_h & = Py \\
P & = I_n - (A')^* A'
\end{align*}
where $y$ is an arbitrary vector in ${\cal F}^n$.
This approach may be preferred in digital computations,
when creating storage for $M$ may be an issue.  However,
this variation would not be preferred if a solution is
sought for additional $b$ vectors, since the ROP step
would have to be repeated.

The second variation (Eqn. \eqref{var2}) transforms $Ax=b$ into $A'x=Mb$.
In terms of augmented matrices, one 
applies the ROP to the augmented matrix $[A|I_n]$ to produce $[A'|M]$,
in which the rows of $A'$ again form a quasi-orthonormal list of vectors.
In this case the solution is
\begin{align*}
x_p & = Gb \\
x_h & = Py \\
P & = I_n - GA \\
G & = (A')^* M
\end{align*}
where $y \in {\cal F}^n$.  Finally, note that in this variation it
is still possible to compute $P$ as $I_n - (A')^* A'$, and to
compute $x_p$ as $(A')^* Mb$.  In other words, one doesn't have to 
actually form $G$ to compute the solution.
Finally, while this case requires storage for $M$, it also allows
one to compute a solution for additional vectors $b$ without 
having to repeat the ROP step.

These differences, while trivial for small example problems,
may become significant when the matrix dimensions are large, and computer
storage space is limited.  Otherwise, choosing one
variation over the other is mainly a matter of convenience.
Pseudocode for these variations is given in Appendix B.

\section*{Penrose Conditions}

It is convenient to classify a generalized inverse ($G$) according to
which Penrose identities it satisifes.  In particular, different 
properties of the solution $Gb$ follow
if certain sets of Penrose identities are satisfied.
The first Penrose identity is $AGA=A$, which is seen
to always be true for our solution.
\begin{align*}
AGA & = A(A')^*MA  \\
 & = M^{-1}MA(A')^*A' \\
 & = M^{-1}A'(A')^*A' \\
 & = M^{-1} {\cal I}_{A'} A'  \\
 & = M^{-1}A' \\
 & = A 
\end{align*}
Likewise, the second identity $GAG=G$ is always true.
\begin{align*}
GAG & = [(A')^*M]A[(A')^*M] \\
    & = (A')^*A'(A')^*M \\
    & = (A')^* {\cal I}_{A'} M \\
    & = (A')^*M \\
	  & = G
\end{align*}
The third Penrose identity ($AG = (AG)^*$) is seen to
be problematic
\begin{align*}
AG & = A(A')^* M \\
   & = M^{-1}MA(A')^*M \\
   & = M^{-1}A' (A')^*M \\
   & = M^{-1} {\cal I}_{A'} M
\end{align*}
If $A$ is of full row rank, then ${\cal I}_{A'}$ equals $I_m$, 
and $AG$ becomes $I_m$; the identity is satisfied.
allowing the identity to be satisfied.  However, 
if $A$ is \emph{not} of full row rank, then this 
identity is not true in general.  Finally, the 
fourth Penrose identity $GA = (GA)^*$ is seen to be true:
\begin{align*}
GA & = (A')^* MA \\
   & = (A')^*A' \\
   & = ((A')^* A')^* \\
   & = ((A')^* MA)^* \\
   & = (GA)^*
\end{align*}
In summary, the generalized inverse obtained by this
algorithm is at least a $\{124\}$-inverse.  If it's 
additionally true that $A$ is of full row rank, then 
it becomes a $\{1234\}$-inverse, a.k.a. a Moore-Penrose inverse.
Recall that generalized inverses which are at least
of type $\{14\}$ (which is the case here) yield
minimum norm solutions $Gb$.

\section*{Online Capabilities}

To the author's knowledge, algorithms for solving $Ax=b$ 
require that all data (i.e., $A$, $b$)
be available at the outset before the linear solver can begin.
This algorithm is different: it can do the calculation as
the data arrives (so long as it arrives in a certain manner).
This style of computation is referred to as an \emph{online}
algorthm.  Examples of when it is 
useful is in cases where it takes a large amount of time to 
compute (or acquire) all the entries in $A$ and $b$.
This approach reduces the overall computation time.  In addition,
it will be shown that the updates to $x$ are mutually
orthogonal; this means that the estimation of $\|x\|$ is
monotonically non-decreasing throughout the computation.

Assume that the rows of $[A|b]$ become available one at a time,
and for simplicty let the i-th row of $[A|b]$ be the i-th to arrive.
Note that when the ROP is based on the classical Gram-Schmidt 
(CGS) \cite{golub,watkins} procedure, the 1st through i-th
rows of $A$, $M$ and $b$ will no longer change following the i-th
step of the algorithm.  Since those rows are done changing at these
points, they become available to be used in a computation of $x_p$.
The next step is to use the column-row expansion \cite{carlson}
on the product $(A')^* b'$ to rewrite $x_p$ as
\begin{equation*}
x_p = (A')^* b' = \sum_{i=1}^m x_p^{(i)}
\end{equation*}
where
\begin{equation*}
x_p^{(i)} = \text{Col}_i[(A')^*] ~ b_i'
\end{equation*}
and Col$_i$ signifies the i-th column.
Following the i-th step in the algorithm, the i-th term on the right-hand side
(ie, $x_p$) may be computed.  Thus the solution is accrued just
by adding $x_p^{(1)} + x_p^{(2)} + \cdots$.
Furthermore, the updates to $x_p$ are mutually orthogonal, i.e.,
$<x_p^{(i)}, x_p^{(k)} > = 0$ for $i \neq k$.  This approach
based on the first variation will be the basis of the following
online computation.  An illustration of this technique
on a numerical example is in Appendix D.

The same basic approach can also be taken for the second
variation of the algorithm, except now the column-row expansion
is used to rewrite the generalized inverse as
\begin{equation*}
G = (A')^* M = \sum_{i=1}^m G^{(i)}
\end{equation*}
where
\begin{equation*}
G^{(i)} = \text{Col}_i[(A')^*] ~ \text{Row}_i[M]
\end{equation*}
Note that $G^{(i)}$ may be computed following the i-th step of the ROP. 
Following the computation of $G$, the particular solution
is easily computed from it.

\section*{Final Remarks}

What is immediately noteworthy about the algorithm\cite{license}
is that it doesn't use an elimination
or partitioning strategy; in particular, there was
no solving of a triangular system of equations.
The algorithm is similar to those based on matrix decompositions \cite{stewart:top10},
in which $A$ is written as a product of matrices, and then
reinserted into $Ax=b$.  In the algorithm presented here, the factorization was done
implicitly, by operating on $A$ while it was in the context
of the equation $Ax=b$.  Borrowing a term from metallurgy,
this type of factorization might be called an "in situ factorization".
Also, keep in mind that it is required that the equations
represented by $Ax=b$ be consistent.  (Although, if they are inconsistent,
a simple modification to the pseudocode makes it easy to discover that.)

The new algorithm is perhaps most similar to a version of the 
QR algorithm, in which the GS procedure operates on
the rows of $A$.  However, in that approach, one still has to
solve a triangular system of equations.  This destroys the
possibility of easily computing a generalized inverse
or a null space projection operator, as well as formulating
the solution in an \emph{online} mode.

The solution found by the new method is always of minimum norm.
When $A$ additionally has full row rank, the solution is a 
least-square solution.  These properties follow from the 
generalized inverse, which is type $\{124\}$ generally, and type $\{1234\}$
when $A$ has full row rank.  (Recall that if a generalized
inverse $G$ is at least type $\{14\}$ the solution $Gb$ has minimum
norm, and if it is at least type $\{13\}$ the solution $Gb$ is a
least squares solution.)  Finally, it's pointed out that it
isn't necessary to explicitly form $G$; the first variation side-steps
that computation.

The orthonormalization procedure named "ROP" can be thought
of as just the GS procedure, except that the resulting zero-norm row vectors
are retained in the end result.  Zero vectors are allowed to
persist in $A$ only because it's easier to leave them there.
They could also be removed; in that case $A$ and $b$ would have
to be re-sized.  An extension of the ROP \cite{zimmer,license}
takes into accout numerical precision and declares the
norm of a vector to be zero if it is less than a small
number $\epsilon$; the size of $\epsilon$ is related to
the machine precision used to implement the algorithm \cite{noble}.

Although the new algorithm was cast to solve $Ax=b$, it can
easily solve \cite{zimmer} its matrix generalization: $AX=B$, where $X$ is n-by-p
and $B$ is m-by-p (and $p \geq 1$).  In that case the ROP proceeds as before,
and the particular and homogeneous parts of the solution 
$X = X_p + X_h$ are
\begin{align*}
 X_p & = G B \\
 X_h & = P Y 
\end{align*}
where $G$ and $P$ are the same as before, and $Y \in {\cal F}^{n \times p}$.
This also admits an online formulation.

\section*{Acknowledgements}
The author dedicates this paper to Robert E. Zimmer for his 
kind support and encouragement.  In addition, the author thanks 
Daniel Grayson for his thoughts and suggestions.

\section*{APPENDIX A: pseudocode}

The following pseudocode, written in the style of the
C programming language, illustrates the workings of the
algorithm.  It shows that while it was previously expedient
to emphasize the role of the matrices $M_s$, it's not
necessary to explicitly form them.  Separate pseudocode is
presented for each of the variations, to aid exposition.

Two accomodations are made for machine precision,
should this be implemented on a computer.  The first
is the variable 'eps', which might be set to some multiple
of the machine precision (cf. "$\epsilon$-rank" \cite{noble}).
(For the examples herein, 'eps' is zero.)
Also, it might be expedient to take further action on a row
that has a norm less than 'eps'; this would be done in the code 
where the comment "zero-norm option" appears.  (However, this
option is not used in the examples herein.)
Finally, the notation $\text{Row}_i[P]$ indicates the i-th row of a matrix $P$.
The bars $||~||$ indicate a Euclidean norm, 
and the angle brackets $<,>$ indicate an inner product.  
Since exact arithmetic is assumed for the examples
and the further discussion in this paper, take 'eps' to be zero.

In the first variation the algorithm begins with the input 
data $A$ and $b$.  Row operations are done on the augmented 
matrix $[A|b]$, transforming it into $[A'|b']$.  The rows of
$A'$ subsequently form a quasi-orthonormal set. This version
of the ROP is based on the modified Gram-Schmidt (MGS) \cite{golub,watkins}
procedure.
\begin{verbatim}
//first variation
for (i = 1 to m){
   mag = || Row_i[A] ||
   if( mag > eps){
      //normalization
      Row_i[A] = Row_i[A] / mag
      b_i = b_i / mag

      //orthogonalization
      for (k = i+1 to m){
         prod = < Row_k[A], Row_i[A] >
         Row_k[A] = Row_k[A] - Row_i[A] * prod
         b_k = b_k - b_i * prod
      }
   }else{
      //implement a "zero-norm option"  
   }
}
\end{verbatim}
Following this the various solution features (i.e., $x_p$, $G$, $P$, ...)
are computed.

In the second variation, the ROP begins with the input
of the coefficient matrix $A$ and an m-by-m matrix $M$
which is initialized as an identity matrix.
In the pseudocode, the entries for $A$ and $M=I_m$ will be written
over; the result at the end will be identified as $A'$
and $M=(\cdots M_2 M_1)$, respectively.
The ROP proceeds by doing row operations on the augmented 
matrix $[A|I_m]$, transforming it into $[A'|M]$, such that 
the rows of $A'$ form a quasi-orthonormal set.
This version of the ROP is also based on the MGS procedure.
\begin{verbatim}
//second variation
for (i = 1 to m){
   mag = || Row_i[A] ||
   if( mag > eps){
      //normalization
      Row_i[A] = Row_i[A] / mag
      Row_i[M] = Row_i[M] / mag

      //orthogonalization
      for (k = i+1 to m){
         prod = < Row_k[A], Row_i[A] >
         Row_k[A] = Row_k[A] - Row_i[A] * prod
         Row_k[M] = Row_k[M] - Row_i[M] * prod
      }
   }else{
      //implement a "zero-norm option"  
   }
}
\end{verbatim}
Following this the various solution features are computed.

\section*{APPENDIX B: example of 1st variation}

In this section the pseudocode for the first variation
is used to compute the solution.  The input data are
\begin{equation*}
A=\begin{bmatrix}
 0 & -3i  &  0\\
2i &   1  & -1 \\
4i & 2-3i & -2 \\
\end{bmatrix}
,\quad\quad
b=\begin{pmatrix}
1 \\
2i \\
1+4i \\
\end{pmatrix}
\end{equation*}
Note that $A$ has rank 2.

\subsection*{Step 1}

The $i=1$ case in the loop in the pseudocode involves the normalization
of the first row.  The associated row operation is
\begin{equation*}
\text{Row}_1 \leftarrow (\frac{1}{3}) \text{Row}_1 
\end{equation*}
Following this operation, the intermediate values for $A'$ and $b'$ are
\begin{equation*}
[A|b]= \begin{bmatrix}
 0  &   -i  &   0  &  | & \frac{1}{3}  \\
2i  &    1  &  -1  &  | &      2i      \\
4i  &  2-3i &  -2  &  | &     1+4i
\end{bmatrix}
\end{equation*}

\subsection*{Step 2}

This step is the orthogonalization between the first and second
rows of $A$; it occurs when $i=1$ and $k=2$ in the pseudocode.  The associated
row operation is
\begin{equation*}
\text{Row}_2 \leftarrow \text{Row}_2  - (i) \text{Row}_1
\end{equation*}
This leads to the following intermediate values for $A'$ and $b'$
\begin{equation*}
[A|b]= \begin{bmatrix}
 0  &  -i   &  0  &  | &   \frac{1}{3}  \\
2i  &   0   & -1  &  | &   \frac{5}{3}i \\
4i  & 2-3i  & -2  &  | &       1+4i
\end{bmatrix}
\end{equation*}

\subsection*{Step 3}

This orthogonalization step is between the first and third
rows, and occurs when i=1 and k=3.  The associated row operation is
\begin{equation*}
\text{Row}_3 \leftarrow \text{Row}_3 - (3+2i) \text{Row}_1
\end{equation*}
At this point the intermediate values for $A'$ and $b'$ are
\begin{equation*}
[A|b]= \begin{bmatrix}
 0 &  -i  &  0  &  | &   \frac{1}{3}  \\
2i &   0  & -1  &  | &   \frac{5}{3}i \\
4i &   0  & -2  &  | &   \frac{10}{3}i
\end{bmatrix}
\end{equation*}

\subsection*{Step 4}

The normalization step for i=2 is for the second row.  The associated
row operation is
\begin{equation*}
\text{Row}_2 \leftarrow  (\frac{1}{\sqrt{5}} ) \text{Row}_2
\end{equation*}
Following this, $A'$ and $b'$ take on the intermediate values
\begin{equation*}
[A|b]= \begin{bmatrix}
        0           &  -i  &          0           &  | &       \frac{1}{3}     \\
\frac{2}{\sqrt{5}}i &   0  & -\frac{1}{\sqrt{5}}  &  | &   \frac{\sqrt{5}}{3}i  \\
        4i          &   0  &         -2           &  | &       \frac{10}{3}i
\end{bmatrix}
\end{equation*}

\subsection*{Step 5}

For i=2 and k=3, the third row is orthogonalized with respect to
the second.  The row operation associated with this is
\begin{equation*}
\text{Row}_3 \leftarrow \text{Row}_3 -  (2\sqrt{5}) \text{Row}_2
\end{equation*}
The intermediate values for $A'$ and $b'$ are now
\begin{equation*}
[A|b]= \begin{bmatrix}
        0           &  -i  &          0           &  | &       \frac{1}{3}     \\
\frac{2}{\sqrt{5}}i &   0  & -\frac{1}{\sqrt{5}}  &  | &   \frac{\sqrt{5}}{3}i  \\
        0           &   0  &          0           &  | &            0
\end{bmatrix}
\end{equation*}
Noteworthy is that the third row of the augmented matrix
has become zero.  Had the third entry in $b'$ been nonzero the equations
would be inconsistent.

\subsection*{Final steps}

According to the pseudocode, the final step should be the
normalization of the third row (for i=3).  However, since
the third row has zero norm, this step is skipped (according
to the if-statement in the pseudocode, since 'eps' is zero).  It follows that what 
were identified as intermediate values for $A'$ and $b'$ in
the fifth step are in fact final values.  The particular solution
may now be formed as
\begin{equation*}
x_p= (A')^* b' = \begin{bmatrix}
 0  &  -\frac{2}{\sqrt{5}}i  &   0   \\
 i  &           0            &   0   \\
 0  &  -\frac{1}{\sqrt{5}}   &   0
\end{bmatrix}
\begin{pmatrix}
\frac{1}{3} \\
\frac{\sqrt{5}}{3}i \\
0 \\
\end{pmatrix}
=\frac{1}{3}\begin{pmatrix}
2  \\
i  \\
-i  
\end{pmatrix}
\end{equation*}
Also, the null space projection operator is
\begin{equation*}
P = I_3 - (A')^*A' = \frac{1}{5}\begin{bmatrix}
 1 & 0 & -2i \\
 0 & 0 &  0  \\
2i & 0 &  4 
\end{bmatrix}
\end{equation*}
As noted earlier, the homogeneous solution is formed as
$x_h=Py$, for arbitrary $y$.  Setting $y=(y_1,y_2,y_3)^T$,
where each entry is arbitrary, the result is
\begin{equation*}
x_h = \frac{1}{5}(y_1 -2i y_3) \begin{pmatrix}
1  \\
0  \\
2i  
\end{pmatrix}
= \alpha \begin{pmatrix}
1  \\
0  \\
2i  
\end{pmatrix}
\end{equation*}
where $\alpha$ is an arbitrary element in ${\cal C}$.
Noteworthy is that the nullity of $A$ is one, which corresponds
to the above parametrization of $x_h$ requiring only one vector.
Also note that it was not necessary to explicitly construct
the $M_s$ that were used in the derivation of the algorithm.

\section*{APPENDIX C: example of 2nd variation}

This is the same example as before, the same steps are
done, except now $M$ is changed instead of $b$.  The row operations
Have the same steps and the same row operation as in the 
example for the first variation.  Hence, all that will be shown
are the intermediate vales for the augmented matrix $[A|M]$.

Step 1
\begin{equation*}
[A|M]= \begin{bmatrix}
 0  &   -i  &   0  &  | &   \frac{1}{3}  &  0  &  0   \\
2i  &    1  &  -1  &  | &        0       &  1  &  0   \\
4i  &  2-3i &  -2  &  | &        0       &  0  &  1
\end{bmatrix}
\end{equation*}

Step 2
\begin{equation*}
[A|M]= \begin{bmatrix}
 0  &  -i   &  0  &  | &   \frac{1}{3}   &  0  &  0   \\
2i  &   0   & -1  &  | &  -\frac{1}{3}i  &  1  &  0   \\
4i  & 2-3i  & -2  &  | &        0        &  0  &  1
\end{bmatrix}
\end{equation*}

Step 3
\begin{equation*}
[A|M]= \begin{bmatrix}
 0 &  -i  &  0  &  | &   \frac{1}{3}     &  0  &  0   \\
2i &   0  & -1  &  | &  -\frac{1}{3}i    &  1  &  0   \\
4i &   0  & -2  &  | &  -1-\frac{2}{3}i  &  0  &  1
\end{bmatrix}
\end{equation*}

Step 4
\begin{equation*}
[A|M]= \begin{bmatrix}
   0    &  -i  &    0  &  |  &   \frac{1}{3}   &  0  &  0  \\
\frac{2}{\sqrt{5}}i & 0 & -\frac{1}{\sqrt{5}} & | 
& -\frac{\sqrt{5}}{15}i & \frac{\sqrt{5}}{5} & 0 \\
  4i    &  0  &  -2  &  |  &  -1-\frac{2}{3}i  &  0  &  1
\end{bmatrix}
\end{equation*}

Step 5
\begin{equation*}
[A|M]= \begin{bmatrix}
   0    &  -i  &    0  &  |  &   \frac{1}{3}   &  0  &  0  \\
\frac{2}{\sqrt{5}}i & 0 & -\frac{1}{\sqrt{5}} & | 
& -\frac{\sqrt{5}}{15}i & \frac{\sqrt{5}}{5} & 0 \\
  0  &  0  &  0  &  |  &   -1  &  -2  &  1
\end{bmatrix}
\end{equation*}
The final $A'$ and $M$ matrices are those found above, in step 5.
The generalized inverse $G$ is found to be
\begin{equation*}
G = (A')^* M = \frac{1}{15}\begin{bmatrix}
-2  &  -6i  &  0  \\
5i  &   0   &  0  \\
i   &  -3   &  0 
\end{bmatrix}
\end{equation*}
Besides permitting a computation of the particular solution
via $x_p = Gb$, it also allows an easy computation of the null
space projection operator ($P = I_n - GA$), which is the same
as before.  In short, the second variation can compute everything
that the first variation did; the difference is that it can also
compute $G$.  Finally, note that it was not necessary to explicitly
construct the $M_s$ that were used in the derivation of the algorithm.

\subsection*{Additional remarks}

To complete the example, the matrices $M_s$ ($s=1,2,3,4,5$)
are shown for the row operations given in the above examples.
($M_s$ corresponds to the row operation in the s-th step.)
While it was not necessary to actually construct and use
these matrices, they may provide an aid to understanding for
the reader.
\begin{equation*}
M_1 = \begin{bmatrix}
\frac{1}{3}  &  0  &  0  \\
 0    &  1  &  0  \\
 0    &  0  &  1  
\end{bmatrix}
,\quad
M_2 = \begin{bmatrix}
  1  &  0  &  0  \\
 -i  &  1  &  0  \\
  0  &  0  &  1  
\end{bmatrix}
,\quad
M_3 = \begin{bmatrix}
     1    &  0  &  0  \\
     0    &  1  &  0  \\
 -(3+2i)  &  0  &  1  
\end{bmatrix},
\end{equation*}
\begin{equation*}
M_4 = \begin{bmatrix}
 1  &      0       &  0  \\
 0  &  \frac{\sqrt{5}}{5}  &  0  \\
 0  &      0       &  1  
\end{bmatrix}
,\quad
M_5 = \begin{bmatrix}
 1  &      0       &  0  \\
 0  &      1       &  0  \\
 0  &  -2\sqrt{5}  &  1  
\end{bmatrix}
\end{equation*}
It is left as an exercise to verify that $M_5 \cdots M_1$
equals the $M$ that was computed in the second variation above.

\section*{APPENDIX D: online case}

In the following example, the same data is used to illustrate
the online case for the first variation.  In the expressions below,
the "input data" for the rows of $A$ and $b$ are from newly
acquired data.
(Double-hyphens in a matrix mean that no data has been entered
there yet.)  Also, the "intermediate results" are how $A$, $b$,
and $x_p$ appear following the $i$-th update.  After $i=3$,
the updates are complete, and the last $A$ and $b$ may be
identified with $A'$ and $b'$, respectively.  Also, note how the
norm of $x_p$ is non-decreasing with respect to $i$. \\

\noindent $i=1$ -----------------------------------------\\
input data:
\begin{align*}
\text{Row}_1(A) & = (0, -3i, 0) \\
\text{Row}_1(b) & = (1)
\end{align*}
intermediate results:
\begin{equation*}
A=\begin{bmatrix}
0   &  -i  &   0 \\
--  &  --  &  -- \\
--  &  --  &  --
\end{bmatrix},\quad
b=\begin{bmatrix}
\frac{1}{3} \\
--  \\
--   
\end{bmatrix}
,\quad
x_p^{(1)}=\begin{bmatrix}
0   \\
\frac{1}{3}i \\
0   
\end{bmatrix}
\end{equation*}
\\

\noindent $i=2$ -----------------------------------------\\
input data:
\begin{align*}
\text{Row}_2(A) & = (2i, 1, -1) \\
\text{Row}_2(b) & = (2i)
\end{align*}
intermediate results:
\begin{equation*}
A=\begin{bmatrix}
0            &  -i  &   0  \\
\frac{2\sqrt{5}}{5}i  &   0  &  -\frac{\sqrt{5}}{5}  \\
--           &  --  &  --
\end{bmatrix},\quad
b=\begin{bmatrix}
\frac{1}{3}   \\
\frac{\sqrt{5}}{3}i  \\
--            
\end{bmatrix}
,\quad
x_p^{(2)}=\begin{bmatrix}
\frac{2}{3}  \\
     0       \\
-\frac{1}{3}i 
\end{bmatrix}
\end{equation*}
\\

\newpage

\noindent $i=3$ -----------------------------------------\\
input data:
\begin{align*}
\text{Row}_3(A) & = (4i, 2-3i, -2) \\
\text{Row}_3(b) & = (1+4i)
\end{align*}
intermediate results:
\begin{equation*}
A=\begin{bmatrix}
0            &  -i  &   0  \\
\frac{2\sqrt{5}}{5}i  &   0  &  -\frac{\sqrt{5}}{5}  \\
0           &  0  &  0
\end{bmatrix},\quad
b=\begin{bmatrix}
\frac{1}{3}   \\
\frac{\sqrt{5}}{3}i  \\
0            
\end{bmatrix}
,\quad
x_p^{(3)} = \begin{bmatrix}
0 \\
0 \\
0
\end{bmatrix}
\end{equation*}
The particular solution follows from adding all the updates, giving
\begin{equation*}
x_p =  x_p^{(1)} + x_p^{(2)} + x_p^{(3)} =  \frac{1}{3} \begin{bmatrix}
2 \\
i \\
-i
\end{bmatrix}
\end{equation*}
which is the same as found earlier.
As a final point, note that it is trivial to repeat the (CGS)
orthonormalization step for each $i$ during this online computation.
Doing so would increase the accuracy of the solution\cite{zimmer}.

\noindent M.~F. Zimmer\\ zim@neomath.net


\begin{thebibliography}{10}
\bibitem{license}
This algorithm, when it appears in certain applications,
has patent protection.  However, that has no bearing on its
non-commercial use in education or research. \\ \\
\end{thebibliography}
\end{document}